\documentstyle[11pt]{article}
\def\Bbb R{{\rm \bf R}}
\def\proclaim#1{\vskip2mm{\bf #1}\em}
\def\endproclaim{\em \vskip2mm}
\def\tag#1{\eqno(#1)}
\def\gathered{\begin{array}{c}}
\def\endgathered{\end{array}}
\def\text{\mbox}

\begin{document}

\title {{\bf The enclosure method for inverse obstacle scattering over a finite time interval: VI. Using shell-type initial data}}
\author{Masaru IKEHATA\footnote{
Laboratory of Mathematics,
Graduate School of Engineering,
Hiroshima University, Higashihiroshima 739-8527, JAPAN}}
%\date{}
\maketitle

\begin{abstract}
A simple idea of finding a domain that encloses an unknown discontinuity embedded in a body
is introduced by considering an inverse boundary value problem for the heat equation.
The idea gives a design of a special heat flux on the surface of the body such that from  
the corresponding temperature field on the surface one can extract the smallest radius of 
the sphere centered at an arbitrary given point in the whole space and enclosing unknown inclusions.  
Unlike before, the designed flux is free from a large parameter.  An application of the idea
to a coupled system of the elastic wave and heat equations are also given.

%{\bf ComputationEnclosureVI.tex}

\noindent
AMS: 35R30, 35R05, 35K05, 35J05, 35L05, 80A23, 74J25, 35Q79, 74F05, 35B40


\end{abstract}

%\tableofcontents
\section{Introduction}
We are interested in seeking an analytical method for inverse obstacle problems in the finite time domain, 
that is, extracting information about the geometry
of unknown discontinuities such as inclusions, cavities and cracks inside a
body by using the observed data on the surface of the body over a finite time interval.  
Here the words ``analytical method'' mean that it fully makes use of the governing equation of 
the observation data in showing the mathematical validity of the method and is based on an {\it extraction 
formula} not just like a conventional optimization method.

The {\it time domain enclosure method} is one of analytical methods and goes back to \cite{I0}.
In \cite{EIV} the method has been applied to an inverse obstacle problem for the wave equation in a three-dimensional bounded domain.
Therein a method using a single set of the input and induced out put data observed on the boundary of the domain over a finite time interval
has been introduced.  The method yields an extraction formula of the distance of an arbitrary fixed point outside the domain to unknown obstacles
inside the domain.

Recently, in \cite{EV} the author introduced a way of combining the method in \cite{EIV} and 
the {\it time-reversal invariance} of the governing equation.
The idea enables us to extract the smallest radius of the sphere centered at an arbitrary given point in the whole space and enclosing unknown obstacles
from a single set of the input and induced out put data observed on the boundary of the domain over a finite time interval.

However, if the governing equation does not have the time-reversal invariance, then one can not apply directly the idea developed in \cite{EV}
to any inverse obstacle problem governed by the equation.   
The heat equation is a typical and an important example.

However, in \cite{EH} we introduced an auxiliary equation
which is a wave equation with a large parameter and using the time-reversal invariance of this equation, we found an extraction formula
of the minimum radius of the sphere centered at an arbitrary given point in the whole space
and enclosing all the unknown cavities inside the body.  Since the wave equation therein contains a large parameter, 
the input heat flux also depends on the same parameter.  
This means that we have to prescribe infinitely many heat fluxes to get one information
about the geometry of the unknown cavity.

In this paper, we introduce a new and extremely simple idea 
which works also for inverse obstacle problems governed by equations without time reversal invariance.  
The idea yields the same information as above by using a single set of input and output data on the surface of the body.
To make the essential difference from the idea in \cite{EV, EH} clear 
and show the applicability to various inverse obstacle problems,
we consider two inverse obstacle problems governed by the heat equation and a coupled system of the elastic wave and heat equations
appearing in the linear theory of thermoelasticity.

\section{Idea}
In this section we explain the idea by considering a typical and important inverse obstacle problem governed by the heat equation.

Let $\Omega$ be a bounded domain of $\Bbb R^3$ with $C^2$-boundary.
Given $f=f(x,t)$ with $x\in\partial\Omega$ and $t\in\,]0,\,T[$ let $u=u_f(x,t)$ with $x\in\Omega$ and $t\in\,]0,\,T[$ solve
$$\left\{
\begin{array}{lll}
(\partial_t-\nabla\cdot\gamma\nabla)u=0, & x\in\Omega, & 0<t<T,\\
\\
\displaystyle
u(x,0)=0, & x\in\Omega, &
\\
\\
\displaystyle
\gamma\nabla u\cdot\nu=f(x,t), & x\in\partial\Omega, & 0<t<T,
\end{array}
\right.
$$
where $\gamma=\gamma(x)$ with $x\in\Omega$ belongs to $L^{\infty}(\Omega)$ and satisfies $\gamma(x)\ge C$ a.e. $x\in\Omega$
for a positive constant $C$; $\nu$ denotes the unit outward normal.

We assume that $\gamma$ takes the form
$$\displaystyle
\gamma(x)=
\left\{
\begin{array}{ll}
\displaystyle
1, & x\in\Omega\setminus D,\\
\\
\displaystyle
1+h(x), & x\in D,
\end{array}
\right.
$$
where $h\in L^{\infty}(D)$ and the set $D$ is an open subset of $\Omega$ with Lipschitz boundary such that $\overline D\subset\Omega$.

Here, we impose two conditions on $h(x)$:

(A.I) There exists $h_{-}>0$ such that $h(x)\le -h_{-}$ a.e. $x\in D$.

(A.II) There exists $h_{+}>0$ such that $h(x)\ge h_{+}$ a.e. $x\in D$.

\noindent
The set $D$ is a mathematical model of unknown inclusions where the conductivity has a negative/positive jump described as (A.I)/(A.II).

We consider the following problem.

$\quad$

{\bf\noindent Problem.}  Assume that both $h$ and $D$ are unknown and that $h$ satisfies (A.I) or (A.II).  
Fix $0<T<\infty$.  
Find a suitable heat flux $f$ in such a way that the temperature field $u_f(x,t)$, $(x,t)\in\partial\Omega\times\,]0,\,T[$
yields some information about the geometry of $D$.

$\quad$

\noindent
In \cite{IK2} we have already shown that if $f$ has a positive lower bound over $\partial\Omega\times\,]0,\,T[$ with
an additional condition, then
one can extract the distance $\text{dist}\,(D,\partial\Omega)=\inf_{x\in D, y\in\partial\Omega}\,\vert x-y\vert$
from $u_f$ observed on $\partial\Omega\times\,]0,\,T[$.  However, this information is too rough, in particular, in the case when 
a piece of $\partial D$ is located near $\partial\Omega$.  Besides, changing the flux $f$ does not yield any effect on the obtained result.

In what follows, we denote by $B_r(y)$ the open ball centered at the point $y\in\Bbb R^3$ with radius $r$;
$\chi_M$ the characteristic function of the set $M\subset\Bbb R^3$.

The construction of the heat flux and the corresponding indicator function in the enclosure method is as follows.

\noindent
(1)  Let $p\in\Bbb R^3$ be an arbitrary point.
Choose $R_1$ in such a way that $\Omega\subset B_{R_1}(p)$
and let $R_2>R_1$.

\noindent
(2)  Solve
$$\left\{
\begin{array}{ll}
\displaystyle
(\partial_t-\Delta)v=0, & (x,t)\in\Bbb R^3\times\,]0,\,\infty[,\\
\\
\displaystyle
v(x,0)=(R_2-\vert x-p\vert)(R_1-\vert x-p\vert)\chi_{B_{R_2}(p)\setminus B_{R_1}(p)}(x), & x\in\Bbb R^3.
\end{array}\right.
$$
Note that the support of the initial data $v(x,0)$ is given by the set $\overline{B_{R_2}(p)}\setminus B_{R_1}(p)$
which is in the form of a spherical shell.

\noindent
(3)  Define the input heat flux $f$ by the equation
$$\begin{array}{lll}
\displaystyle
f(x,t)=\nabla v(x,t)\cdot\nu(x), & x\in\partial\Omega, & 0<t<T,
\end{array}
$$
and the indicator function
$$\begin{array}{ll}
\displaystyle
I(\tau;R_1,R_2,p)=\int_{\partial\Omega}(w_0-w)\nabla w_0\cdot\nu\,dS, & \tau>0,
\end{array}
$$
where
$$\left\{
\begin{array}{ll}
\displaystyle
w_0(x;\tau)=\int_0^Te^{-\tau t}v(x,t)\,dt, & x\in\Bbb R^3
\\
\\
\displaystyle
w(x;\tau)=\int_0^Te^{-\tau t}u_f(x,t)\,dt, & x\in\Omega.
\end{array}
\right.
$$
Note that the solution $v$ has the form
$$\begin{array}{lll}
\displaystyle
v(x,t)=\frac{1}{(\sqrt{2\pi t}\,)^3}
\int_{B_{R_2}(p)\setminus B_{R_1}(p)}e^{-\frac{\vert x-y\vert^2}{4t}}(R_2-\vert y-p\vert)(R_1-\vert y-p\vert)\,dy,
&
x\in\Bbb R^3,
& t>0,
\end{array}
$$
and thus flux $f$ is given explicitly.

\noindent
The information about the geometry of the unknown inclusion considered in this paper is the quantity defined by
$$\displaystyle
R_D(p)=\sup_{x\in D}\,\vert x-p\vert.
$$
This is the minimum radius of the sphere cenetered at point $p$ and enclosing $D$.  Thus we have
$$\displaystyle
D\subset B_{R_D(p)}(p).
$$
See Figure 1 for an illustration of the ball $B_{R_D(p)}(p)$.

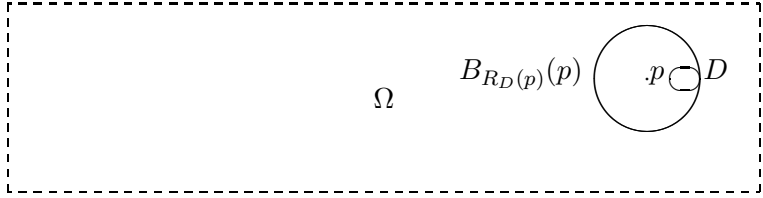
\begin{figure}
\setlength\unitlength{0.5truecm}
  \begin{picture}(10,10)(-3,-1)
  \put(3,0){\dashbox{0.2}(20,5){$\Omega$}}
  \put(21,3){\oval(0.8,0.6)}
  \put(15,3){$B_{R_D(p)}(p)$}
  \put(20,3){\circle{5}}
  \put(20,3){\circle{5}}
  \put(21.5,3){$D$}
  \put(20,3){\circle*{0.1}{$p$}}
  \end{picture}
 \caption{\label{fig:1}  An illustration of ball $B_{R_D(p)}(p)$}
 \end{figure}

The following theorem is the main result of this paper.
\proclaim{\noindent Theorem 2.1.}  Let $p$ be an arbitrary point in $\Bbb R^3$.
Then, from the asymptotic behaviour of indicator function $I(\tau;R_1,R_2,p)$ as $\tau\rightarrow\infty$ one has the value of 
$R_D(p)$.
\endproclaim

\noindent
Note that the heat flux $f$ does not contain the parameter $\tau$ which is the independent variable
of the indicator function.
Only in the process of the data, in other words, in computing the indicator function 
after having the observation data one makes variable $\tau$ large.  Thus a single set of input and output data yields one information.  
This gives a complete answer to the question raised in \cite{EH}.

Here we give a short proof of Theorem 2.1 focused on the main point of the idea.

We see that the function $w_0$ satisfies
$$\begin{array}{ll}
\displaystyle
(\Delta-\tau)w_0+v(x,0)=e^{-\tau T}v(x,T), & x\in\Bbb R^3.
\end{array}
\tag {2.1}
$$
Besides, since we have chosen $R_1$ in such a way that $\Omega\subset B_{R_1}(p)$, initial data
$v(x,0)$ vanishes for $x\in\Omega$.  Thus one gets
$$\begin{array}{ll}
\displaystyle
(\Delta-\tau)w_0=e^{-\tau T}v(x,T), & x\in\Omega.
\end{array}
$$
On the other hand, the function $w$ satisfies
$$\left\{
\begin{array}{ll}
\displaystyle
(\nabla\cdot\gamma\nabla-\tau)w=e^{-\tau T}u_f(x,T), & x\in\Omega,
\\
\\
\displaystyle
\gamma\nabla w(x;\tau)\cdot\nu(x)=\nabla w_0(x;\tau)\cdot\nu(x)=\int_0^Te^{-\tau t}f(x,t)dt, & x\in\partial\Omega.
\end{array}
\right.
$$
Then, from the previous studies \cite{IK2, IHF}, roughly speaking, we see that, as $\tau\rightarrow\infty$
$$\displaystyle
I(\tau;R_1,R_2,p)\sim \mp \int_D\vert\nabla w_0\vert^2\,dx,
$$
where $-/+$ should be chosen according to (A.I)/(A.II).  This means that the asymptotic behaviour of
the indicator function is governed by that of the right-hand side.

From (2.1) together with $w_0\in H^1(\Bbb R^3)$ we see that
$$\displaystyle
w_0(x;\tau)\sim w_{00}(x;\tau)
$$
where
$$\begin{array}{ll}
\displaystyle
w_{00}(x;\tau)
=\frac{1}{4\pi}\int_{\Bbb R^3}\frac{e^{-\sqrt{\tau}\,\vert x-y\vert}}{\vert x-y\vert}\,v(y,0)\,dy, & x\in\Bbb R^3.
\end{array}
$$
Thus the key point is the choice of $v(x,0)$.  
Since we have chosen as
$$\displaystyle
v(x,0)=(R_2-\vert x-p\vert)(R_1-\vert x-p\vert)\chi_{B_{R_2}(p)\setminus B_{R_1}(p)}(x),
$$
the function $w_{00}$ takes the form
$$\displaystyle
w_{00}(x;\tau)=\frac{1}{4\pi}\int_{B_{R_2}(p)\setminus B_{R_1}(p)}\frac{e^{-\sqrt{\tau}\,\vert x-y\vert}}{\vert x-y\vert}\,
(R_2-\vert y-p\vert)(R_1-\vert y-p\vert)
\,dy.
\tag {2.2}
$$
For this we have
\proclaim{\noindent Lemma 2.2.}
We have the expression
$$\begin{array}{ll}
\displaystyle
w_{00}(x;\tau)=H(\tau;R_1,R_2)e^{-\sqrt{\tau}\,R_1}\frac{\sinh\,\sqrt{\tau}\vert x-p\vert}{\vert x-p\vert}, & x\in B_{R_1}(p),
\end{array}
\tag {2.3}
$$
where the coefficient $H(\tau;R_1,R_2)$ is independent of $x$ and satisfies
$$\displaystyle
\lim_{\tau\rightarrow\infty}(\sqrt{\tau})^3\,H(\tau;R_1,R_2)=-R_1(R_2-R_1).
$$
\endproclaim

Note that the function
$$\displaystyle
\Bbb R^3\setminus\{p\}\ni x\mapsto \frac{\sinh\sqrt{\tau}\,\vert x-p\vert}{\vert x-p\vert}
$$
has the unique extension to the whole space as the solution of the modified Helmholtz equation
$(\Delta-\tau)v=0$.  The function above appearing in the right-hand side on (2.3) 
should be the one replaced with this extension.  Throughout this paper, we always make use of this convention.

For each $x\in B_{R_1}(p)$,
$$\displaystyle
\inf_{y\in B_{R_2}(p)\setminus B_{R_1}(p)}\,\vert x-y\vert=R_1-\vert x-p\vert.
$$
This together with (2.2) suggests, by the Laplace method, the asymptotic profile of $w_{00}(x;\tau)$
may take the form $e^{-\sqrt{\tau}\,(R_1-\vert x-p\vert)}$ multiplied by an algebraic power of $\tau$ as $\tau\rightarrow\infty$.
Its precise expression is given by the formula (2.3).

Now admit Lemma 2.2 and move on.  Then, since $\Omega\subset B_{R_1}(p)$, for $x\in D$ we will have, as $\tau\rightarrow\infty$
$$
\displaystyle
\nabla w_{00}(x;\tau)\sim \tau^{-1}e^{-\sqrt{\tau}\, R_1}e^{\tau\vert x-p\vert}.
$$

\noindent
This gives
$$\displaystyle
\int_D\vert\nabla w_{00}(x;\tau)\vert^2 dx\sim e^{-2\sqrt{\tau}\, R_1}e^{2\sqrt{\tau} R_D(p)}.
$$
See Lemma 2.4 in \cite{EV} for the proof together with its exact meaning.  Note that 
we ignore some algebraic growing or decaying factor with respec to $\tau$ in the right-hand side
and the following also adopts this rule.

Therefore one gets, as $\tau\rightarrow\infty$
$$\displaystyle
I(\tau;R_1,R_2,p)\sim\mp e^{-2\sqrt{\tau} R_1}e^{2\sqrt{\tau} R_D(p)}.
$$
This yields the extraction formula of $R_D(p)$:
$$\displaystyle
\lim_{\tau\rightarrow\infty}\frac{1}{\sqrt{\tau}}\,\log\vert I(\tau;R_1,R_2)\vert=2(R_D(p)-R_1).
$$
This gives a {\it quantitative} information about the geometry of $D$.

Moreover, we have
$$\displaystyle
\lim_{\tau\rightarrow\infty}e^{\sqrt{\tau}\,T}I(\tau;R_1,R_2)
=
\left\{
\begin{array}{ll}
\displaystyle
0 & \text{if $T<2(R_1-R_D(p))$,}
\\
\\
\displaystyle
\mp\infty & \text{if $T>2(R_1-R_D(p))$ and (A.I)/(A.II) is satisfied.}
\end{array}
\right.
$$

Thus, there exists a positive number $\tau_0$ such that

$\quad$

$\bullet$ if (A.I) is satisfied, then for all $\tau\ge\tau_0$ we have $I(\tau;R_1,R_2)<0$,

$\quad$

$\bullet$ if (A.II) is satisfied, then for all $\tau\ge\tau_0$ we have $I(\tau;R_1,R_2)>0$.

$\quad$

\noindent
This gives a {\it qualitative} target distinction.

These are the proof and exact statement of Theorem 2.1.  The proof of Lemma 2.2 is given in Section 3.

\subsection{Comparison}

Here we make a comparison of the idea developed in this paper with that of \cite{EH}. 
In \cite{EH}, to generate a suitable heat flux we solve the Cauchy problem for an auxiliary wave equation 
with a parameter.

More precisely, given $p\in\Bbb R^3$ and $\eta>0$, let $v=v(x,t;\tau)$ be the solution of
$$\left\{
\begin{array}{lll}
\displaystyle
(\partial_t^2-\tau\Delta)v=0, & x\in\Bbb R^3, & t>0,\\
\\
\displaystyle
v(x,0;\tau)=0, & x\in\Bbb R^3, &\\
\\
\displaystyle
\partial_tv(x,0;\tau)=(\eta-\vert x-p\vert)\,\chi_{B_{\eta}(p)}, & x\in\Bbb R^3. &
\end{array}
\right.
$$
where $\tau$ is the same positive parameter as above.  Note that the propagation speed of the wave equation
grows to infinity as $\tau\rightarrow\infty$.

Using the solution $v=v(x,t;\tau)$, we prescribe the heat flux $f$ given by
$$\begin{array}{lll}
\displaystyle
f(x,t;\tau)=\nabla v(x,T-t;\tau)\cdot\nu(x),  & x\in\partial\Omega, & 0<t<T.
\end{array}
$$
The points of this formula consist of two parts: 

(i)  prescribing the heat flux coming from a wave equation with a parameter over a finite time interval which goes back to \cite{EH0};

(ii) the {\it time reversal operation} $t\mapsto T-t$ is done for the solution in (i) \cite{EV}.

Comparing with this choice, we see that the heat flux in Theorem 2.1 is naturally constructed and thus, to solve an 
inverse problems for the heat equation we do not make use of any other auxiliary equation in the time domain.
Besides, the heat flux given above depends on a parameter, thus in this sense, 
{\it infinitely many} heat fluxes are prescribed.

Instead of $w_0(x;\tau)$ we employ the function $w^*(x;\tau)$ defined by
$$\begin{array}{ll}
\displaystyle
w^*(x;\tau)=\int_0^T e^{-\tau t}v(x,T-t;\tau)dt, & x\in\Bbb R^3.
\end{array}
$$
Function $w^*$ satisfies
$$\begin{array}
{ll}
\displaystyle
(\Delta-\tau)w^*+\tau^{-1}(\tau v(x,T;\tau)-\partial_tv(x,T;\tau))=\tau^{-1}e^{-\tau T}(\eta-\vert x-p\vert)\chi_{B_{\eta}(p)}(x),
& x\in\Bbb R^3.
\end{array}
$$
Comparing this with (2.1), we see that the term involving functions $v(x,T;\tau)$ and $\partial_tv(x,T;\tau)$ play
the role of the initial data $v(x,0)$ appeared in (2.1).

We can see that
$$\displaystyle
w^*(x;\tau)\sim w^{**}(x;\tau),
$$
where
$$\begin{array}{ll}
\displaystyle
w^{**}(x;\tau)=\frac{1}{4\pi}
\int_{\Bbb R^3}\frac{e^{-\sqrt{\tau}\,\vert x-y\vert}}{\vert x-y\vert}\,\tau^{-1}(\tau v(y,T;\tau)-\partial_tv(y,T;\tau))\,dy, & x\in\Bbb R^3.
\end{array}
$$
Then the Huygens principle tells us that 
$$\displaystyle
v(x,T;\tau)=\partial_t v(x,T;\tau)=0
$$
for $x\in B_{\sqrt{\tau}T-\eta}(p)\cup (\Bbb R^3\setminus B_{\sqrt{\tau}T+\eta}(p))$.
Thus, this becomes
$$\begin{array}{ll}
\displaystyle
w^{**}(x;\tau)=\frac{1}{4\pi}
\int_{B_{\sqrt{\tau}T+\eta}(p)\setminus B_{\sqrt{\tau}T-\eta}(p)}\frac{e^{-\sqrt{\tau}\,\vert x-y\vert}}{\vert x-y\vert}\,\tau^{-1}
(\tau v(y,T;\tau)-\partial_tv(y,T;\tau))\,dy, & x\in\Bbb R^3.
\end{array}
$$
Thus comparing this with (2.2) we will see that function $w_{00}(x;\tau)$ may play the same role of the function $w^{**}(x;\tau)$.
In fact, the profile of $w^{**}$ has been given by
$$\begin{array}{ll}
\displaystyle
w_1^*(x;\tau)=H^*(\tau,T,\eta)\,e^{-\sqrt{\tau}(\sqrt{\tau}\,T-\eta)}\,\frac{\sinh\,\sqrt{\tau}\,\vert x-p\vert}{\vert x-p\vert}, & x\in B_{\sqrt{\tau}T-\eta}(p),
\end{array}
\tag {2.4}
$$
where $\tau$ is large in the sense that
$\Omega\subset B_{\sqrt{\tau}T-\eta}(p)$; the coefficient $H^*(\tau,T,\eta)$ is independent of $x$ and satisfies
$$\displaystyle
\lim_{\tau\rightarrow\infty}(\sqrt{\tau})^3\tau H^*(\tau;T,\eta)=\eta.
$$
However, for establishing (2.4) we need the detailed knowledge of the solution $v(x,t;\tau)$ together with $\partial_tv(x,t;\tau)$
at $t=T$ and near $\partial B_{\sqrt{\tau}T-\eta}(p)$.  See Proposition 3.3 in \cite{EV} whose proof is based on Kirchhoff's formula.
Another proof is based on the representation of the Fourier transform of $v(x,T;\tau)$ together with $\partial_tv(x,T;\tau)$
with respect to $x$.  See the proof of Proposition 3.1 in \cite{EH}.

From these comparison, it is clear that the idea developed in this paper is much simpler than that of
\cite{EH,EV} and can be directly applied to various inverse obstacle problems where the governing equations 
do not have the time-reversal invariance.

\section{Proof of Lemma 2.2}

In this section, for simplicity of description, we replace $\sqrt{\tau}$ with $\tau$.
Let $0<R_1<R_2$ and $x\in B_{R_1}(p)$.

First we compute the volume potential
$$\displaystyle
I(x)=\frac{1}{4\pi}\int_{B_{R_2}(p)\setminus B_{R_1}(p)}
\frac{e^{-\tau\vert x-y\vert}}{\vert x-y\vert}\,(R_2-\vert y-p\vert)(R_1-\vert y-p\vert)\,dy.
$$
One can write
$$\displaystyle
I(x)=I_2(x)-I_1(x),
$$
where
$$\left\{
\begin{array}{l}
\displaystyle
I_1(x)=\frac{1}{4\pi}\int_{B_{R_1}(p)}
\frac{e^{-\tau\vert x-y\vert}}{\vert x-y\vert}\,(R_2-\vert y-p\vert)(R_1-\vert y-p\vert)\,dy,
\\
\\
\displaystyle
I_2(x)=\frac{1}{4\pi}\int_{B_{R_2}(p)}
\frac{e^{-\tau\vert x-y\vert}}{\vert x-y\vert}\,(R_1-\vert y-p\vert)(R_2-\vert y-p\vert)\,dy.
\end{array}
\right.
$$
Thus it suffices to compute $I_j(x)$ for $x\in B_{R_j}(p)$.

Since
$$\displaystyle
(R_2-\vert y-p\vert)(R_1-\vert y-p\vert)
=R_1R_2-(R_1+R_2)\vert y-p\vert+\vert y-p\vert^2,
$$
one has the expression
$$\left\{
\begin{array}{l}
\displaystyle
I_1(x)=R_1R_2 v_0^1(x)-(R_1+R_2)v_1^1(x)+v_2^1(x),
\\
\\
\displaystyle
I_2(x)=R_1R_2 v_0^2(x)-(R_1+R_2)v_1^2(x)+v_2^2(x),
\end{array}
\right.
$$
where
$$\displaystyle
v_j^i(x)=\frac{1}{4\pi}\int_{B_{R_i}(p)}\frac{e^{-\tau\vert x-y\vert}}{\vert x-y\vert}\,\vert y-p\vert^j\,dy.
$$
Therefore we have
$$\displaystyle
I(x)=R_1R_2(v_0^2(x)-v_0^1(x))-(R_1+R_2)(v_1^2(x)-v_1^1(x))+(v_2^2(x)-v_2^1(x)).
$$
The previous computation (see Proposition A.1 in Appendix of \cite{EV}) gives
$$\begin{array}{l}
\displaystyle
\,\,\,\,\,\,
\tau^2(v_0^2(x)-v_0^1(x))
\\
\\
\displaystyle
=\left\{\left(R_1+\frac{1}{\tau}\right)e^{-\tau R_1}-\left(R_2+\frac{1}{\tau}\right)e^{-\tau R_2}\right\}
\frac{\sinh\tau\vert x-p\vert}{\vert x-p\vert}\\
\\
\displaystyle
=\left\{\left(R_1+\frac{1}{\tau}\right)-\left(R_2+\frac{1}{\tau}\right)e^{-\tau (R_2-R_1)}\right\}e^{-\tau R_1}
\frac{\sinh\tau\vert x-p\vert}{\vert x-p\vert}\\
\\
\displaystyle
=\left\{\left(R_1+\frac{1}{\tau}\right)+O(e^{-\tau (R_2-R_1)})\right\}e^{-\tau R_1}
\frac{\sinh\tau\vert x-p\vert}{\vert x-p\vert},
\end{array}
$$
$$\begin{array}{l}
\displaystyle
\,\,\,\,\,\,
\tau^2(v_1^2(x)-v_1^1(x))
\\
\\
\displaystyle
=\left\{\left(R_1^2+\frac{2}{\tau}R_1+\frac{2}{\tau^2}\right)e^{-\tau R_1}
-\left(R_2^2+\frac{2}{\tau}R_2+\frac{2}{\tau^2}\right)e^{-\tau R_2}
\right\}\frac{\sinh\tau\vert x-p\vert}{\vert x-p\vert}
\\
\\
\displaystyle
=\left\{\left(R_1^2+\frac{2}{\tau}R_1+\frac{2}{\tau^2}\right)
+O(e^{-\tau(R_2-R_1)})
\right\}
e^{-\tau R_1}\frac{\sinh\tau\vert x-p\vert}{\vert x-p\vert}
\end{array}
$$
and
$$\begin{array}{l}
\displaystyle
\,\,\,\,\,\,
\tau^2(v_2^2(x)-v_2^1(x))
\\
\\
\displaystyle
=\left\{
\left(R_1^3+\frac{3R_1^2}{\tau}+\frac{6R_1}{\tau^2}+\frac{6}{\tau^3}\right)e^{-\tau R_1}
-\left(R_2^3+\frac{3R_2^2}{\tau}+\frac{6R_2}{\tau^2}+\frac{6}{\tau^3}\right)e^{-\tau R_2}
\right\}\frac{\sinh\tau\vert x-p\vert}{\vert x-p\vert}
\\
\\
\displaystyle
=\left\{
\left(R_1^3+\frac{3R_1^2}{\tau}+\frac{6R_1}{\tau^2}+\frac{6}{\tau^3}\right)
+O(e^{-\tau(R_2-R_1)})
\right\}e^{-\tau R_1}
\frac{\sinh\tau\vert x-p\vert}{\vert x-p\vert}.
\end{array}
$$
Thus we see that, as $\tau\rightarrow\infty$
$$\displaystyle
I(x)\sim \frac{R_1(R_1-R_2)}{\tau^3}e^{-\tau R_1}\frac{\sinh\tau\vert x-p\vert}{\vert x-p\vert}.
$$
This completes the proof of Lemma 2.2.

\section{Application to thermoelasticity}

In \cite{IThermo} the time domain enclosure method developed in \cite{EIV} has been applied to an inverse obstacle problem 
governed by a classical system in the linear theory of thermoelasticity, see \cite{C} for the whole knowledge 
about the system.

It is a coupled system of the elastic wave and heat equations.  Needless to say, the system does not have the
time reversal invariance.  So it is interesting to consider whether the idea explained in the former sections
can be applied to the inverse obstacle problem or not.

Now let us describe the governing system.
In what follows, we assume that $\Omega\setminus\overline D$ is connected.
Let $0<T<\infty$.  Given $f=f(x,t)$ and $\mbox{\boldmath $G$}=\mbox{\boldmath $G$}(x,t)$ with
$(x,t)\in\partial\Omega\times\,]0,\,T[$, let
$\mbox{\boldmath $u$}=
\mbox{\boldmath $u$}_{f,\mbox{\boldmath $G$}}(x,t)$
and $\vartheta=\vartheta_{f, \mbox{\boldmath $G$}}(x,t)$ with $(x,t)\in\,(\Omega\setminus\overline D)\times\,]0,\,T[$ 
denote the solutions of the following initial boundary value problem
$$\displaystyle
\left\{
\begin{array}{lll}
\displaystyle
\rho\partial_t^2\mbox{\boldmath $u$}-\mu\Delta\mbox{\boldmath $u$}
-(\lambda+\mu)\nabla(\nabla\cdot\mbox{\boldmath $u$})-m\nabla\vartheta=\mbox{\boldmath $0$},
& x\in\,\Omega\setminus\overline D, & 0<t<T,\\
\\
\displaystyle
c\partial_t\vartheta-k\Delta\vartheta-m\theta_0\nabla\cdot\partial_t\mbox{\boldmath $u$}=0,
& x\in\,\Omega\setminus\overline D, & 0<t<T,
\\
\\
\displaystyle
\mbox{\boldmath $u$}(x,0)=\mbox{\boldmath $0$}, & x\in\,\Omega\setminus\overline D, &
\\
\\
\displaystyle
\partial_t\mbox{\boldmath $u$}(x,0)=\mbox{\boldmath $0$}, & x\in\,\Omega\setminus\overline D, &
\\
\\
\displaystyle
\vartheta(x,0)=0, & x\in\,\Omega\setminus\overline D, &
\\
\\
\displaystyle
s(\mbox{\boldmath $u$}, \vartheta)\mbox{\boldmath $\nu$}
=\mbox{\boldmath $0$}, & x\in\,\partial D, & 0<t<T,\\
\\
\displaystyle
k\nabla\vartheta\cdot\mbox{\boldmath $\nu$}=0,
& x\in\,\partial D, & 0<t<T,\\
\\
\displaystyle
s(\mbox{\boldmath $u$},\vartheta)
\mbox{\boldmath $\nu$}=\mbox{\boldmath $G$}, & x\in\,\partial\Omega, & 0<t<T,
\\
\\
\displaystyle
-k\nabla\vartheta\cdot\mbox{\boldmath $\nu$}=f, & x\in\,\partial\Omega, & 0<t<T,
\end{array}
\right.
\tag {4.1}
$$
where
$$\
\displaystyle
s(\mbox{\boldmath $u$},\vartheta)
=2\mu\,\text{Sym}\,\nabla\mbox{\boldmath $u$}+
\lambda(\nabla\cdot\mbox{\boldmath $u$})I_3+m\,\vartheta I_3.
$$
The constants
$\theta_0$ and $m$ are the reference temperature and stress-temperature modulus of the body
$\Omega\setminus\overline D$, respectively;
$k$ the conductivity; $\lambda$ and $\mu$ are Lam\'e modulus and shear modulus, respectively;
$\rho$ and $c$ the density and specific heat.   
$\vartheta$ denotes the temperature difference of the absolute temperature from the reference temperature
$\theta_0$;
$\mbox{\boldmath $u$}$ and $s(\mbox{\boldmath $u$}, \vartheta)\mbox{\boldmath $\nu$}$ the displacement vector field
and the surface traction, respectively.

It is assumed that $\rho$, $c$, $\theta_0$ and $k$ are known positive constants, and $m$, $\lambda$ and $\mu$ are known
constants and satisfy $m\not=0$, $\mu>0$ and $3\lambda+2\mu>0$.

The direct problem has been studied in \cite{D} in a more general setting under the assumption that both $\partial\Omega$ and
$\partial D$ are smooth.  In this section we employ this smoothness assumption.

The inverse obstacle problem considered in \cite{IThermo} is the following.

$\quad$

{\bf\noindent Problem.}  Fix a large $T$ (to be determined later) and a single set of the admissible 
pair  $(f, \mbox{\boldmath $G$})$ (to be specified later).
Assume that set $D$ is  unknown.  Extract information about the
location and shape of $D$ from the displacement field $\mbox{\boldmath $u$}_{f,\mbox{\boldmath $G$}}(x,t)$
and temperature difference $\vartheta_{f,\mbox{\boldmath $G$}}(x,t)$ given at
all $x\in\partial\Omega$ and $t\in\,]0,\,T[$.

$\quad$

\noindent
Applying the idea in \cite{EIV} to the system (4.1), we have obtained an extraction formula of the distance of $D$ to an arbitrary fixed point 
outside $\Omega$.  Here we apply the present method which enables us to extract the value of $R_D(p)$ for an arbitrary fixed point $p$ in the whole space.

Let $R_1$ and $R_2$ be the same as those in the preceeding sections.
Given an arbitrary unit vector $\mbox{\boldmath $a$}$, let $\mbox{\boldmath $\Phi$}$ solve
$$\left\{
\begin{array}{lll}
(\rho\partial_t^2-\mu\Delta)\mbox{\boldmath $\Phi$}=\mbox{\boldmath $0$}, & x\in\Bbb R^3, & 0<t<T,\\
\\
\displaystyle
\mbox{\boldmath $\Phi$}(x,0)=\mbox{\boldmath $0$}, & x\in\Bbb R^3, & \\
\\
\displaystyle
\partial_t\mbox{\boldmath $\Phi$}(x,0)=(R_1-\vert x-p\vert)^2(R_2-\vert x-p\vert)^2\chi_{B_{R_2}(p)\setminus B_{R_1}(p)}(x)
\mbox{\boldmath $a$}, & x\in\Bbb R^3. &
\end{array}
\right.
\tag {4.2}
$$
Define
$$\displaystyle
\mbox{\boldmath $v$}_s=\nabla\times\mbox{\boldmath $\Phi$}.
\tag {4.3}
$$
Since we have $\nabla\cdot\mbox{\boldmath $v$}_s=\mbox{\boldmath $0$}$, 
we see that the pair $(\mbox{\boldmath $v$},\Theta)$ given by $(\mbox{\boldmath $v$}_s,0)$ satisfies
$$\displaystyle
\left\{
\begin{array}{lll}
\displaystyle
\rho\partial_t^2\mbox{\boldmath $v$}-\mu\Delta\mbox{\boldmath $v$}
-(\lambda+\mu)\nabla(\nabla\cdot\mbox{\boldmath $v$})-m\nabla\Theta=\mbox{\boldmath $0$},
& x\in\,\Bbb R^3, & 0<t<T,\\
\\
\displaystyle
c\partial_t\Theta-k\Delta\Theta-m\theta_0\nabla\cdot\partial_t\mbox{\boldmath $v$}=0,
& x\in\,\Bbb R^3, & 0<t<T,\\
\\
\displaystyle
\mbox{\boldmath $v$}(x,0)=\mbox{\boldmath $0$}, & x\in\Bbb R^3, &
\\
\\
\displaystyle
\partial_t\mbox{\boldmath $v$}(x,0)=\nabla\times\partial_t\mbox{\boldmath $\Phi$}(x,0), & x\in\Bbb R^3, & 
\\
\\
\displaystyle
\Theta(x,0)=0, & x\in\Bbb R^3. &
\end{array}
\right.
$$
The input data $\mbox{\boldmath $G$}$ and $f$ in (4.1) is given by
$$\left\{
\begin{array}{lll}
\mbox{\boldmath $G$}(x,t)=s(\mbox{\boldmath $v$}_s,0)(x,t)\mbox{\boldmath $\nu$}(x), & x\in\partial\Omega, & 0<t<T,\\
\\
\displaystyle
f(x,t)=0, & x\in\partial\Omega, & 0<t<T.
\end{array}
\right.
\tag {4.4}
$$

We compute the indicator function $I^1(\tau;\mbox{\boldmath $v$}_s,0)$ defined by (1.6) in \cite{IThermo}
for the solution $\mbox{\boldmath $u$}$ of (4.1) with $f$ and $\mbox{\boldmath $G$}$ given by (4.4):
$$\displaystyle
I^1(\tau;\mbox{\boldmath $v$}_s,0)
=\int_{\partial\Omega}s(\mbox{\boldmath $w$}_0,0)\mbox{\boldmath $\nu$}\cdot(\mbox{\boldmath $w$}-\mbox{\boldmath $w$}_0)\,dS,
$$
where
$$\left\{
\begin{array}{ll}
\displaystyle
\mbox{\boldmath $w$}(x;\tau)=\int_0^Te^{-\tau t}\mbox{\boldmath $u$}(x,t)dt, & x\in\Omega\setminus\overline D,\\
\\
\displaystyle
\mbox{\boldmath $w$}_0(x;\tau)=\int_0^Te^{-\tau t}\mbox{\boldmath $v$}_s(x,t)dt, & x\in\Bbb R^3.\\
\end{array}
\right.
$$

\noindent
Then, from the asymptotic behaviour of indicator function $I^1(\tau;\mbox{\boldmath $v$}_s,0)$ as $\tau\rightarrow\infty$,
one obtains the value of $R_D(p)$ at an arbitrary point $p\in\Bbb R^3$.
More precisely, we have

\proclaim{\noindent Theorem 4.1.}  Let $0<R_1<R_2$ and $\Omega\subset B_{R_1}(p)$.

(i)  Let $T$ satisfy
$$\displaystyle
T>\sqrt{\frac{\rho}{\mu}}\,(R_1-R_D(p)+R_{\Omega}(p)-R_D(p)).
\tag {4.5}
$$

Then, there exists a positive number $\tau_0$ such that $I^1(\tau;\mbox{\boldmath $v$}_s,0)>0$ for all $\tau\ge\tau_0$,
and we have
$$\displaystyle
\lim_{\tau\rightarrow\infty}\,\frac{1}{\tau}\log I^1(\tau;\mbox{\boldmath $v$}_s,0)=-2\sqrt{\frac{\rho}{\mu}}(R_1-R_D(p)).
\tag {4.6}
$$

(ii)  We have
$$\displaystyle
\lim_{\tau\rightarrow\infty}e^{\tau T}I^1(\tau;\mbox{\boldmath $v$}_s,0)
=
\left\{
\begin{array}{ll}
\displaystyle
\infty 
&
\displaystyle
\text{if 
$\displaystyle T>2\sqrt{\frac{\rho}{\mu}}\,(R_1-R_D(p))$,}
\\
\\
\displaystyle
0
&
\displaystyle
\text{if 
$\displaystyle T\le 2\sqrt{\frac{\rho}{\mu}}\,(R_1-R_D(p))$.}
\end{array}
\right.
\tag {4.7}
$$

\endproclaim

\noindent
Three remarks are in order.

(a)  The quantity of the right-hand side on (4.5) can be viewed {\it intuitively} as being almost equal to a travel time of the signal
with the speed $\sqrt{\mu/\rho}$ passing through
$$\displaystyle
\partial B_{R_1}(p)\longrightarrow\partial D\longrightarrow\partial\Omega.
$$ 
For example, consider the case when $D=B_{R_D(p)}(p)$ and $\Omega=B_{R_{\Omega}(p)}(p)$.

(b)  Condition (4.5) depends on the size of $R_D(p)$ which should be unknown.
Since 
$$\displaystyle
R_1-R_D(p)+R_{\Omega}(p)-R_D(p)<R_1+R_{\Omega}(p),
$$
we see that if the time $T$ satisfies
$$\displaystyle
T\ge\sqrt{\frac{\rho}{\mu}}(R_1+R_{\Omega}(p)),
$$
then (4.5) is satisfied with the $T$.  Since $R_1$ satisfies $\Omega\subset B_{R_1}(p)$, the minimum choice
of $R_1$ is $R_1=R_{\Omega}(p)$.  Thus, we have the formula (4.6) under the condition
$$\displaystyle
T\ge 2\sqrt{\frac{\rho}{\mu}}R_{\Omega}(p),
$$
which is {\it independent} of unknown $R_D(p)$.

(c) Note that we have
$$\displaystyle
2\sqrt{\frac{\rho}{\mu}}\,(R_1-R_D(p))\ge \sqrt{\frac{\rho}{\mu}}\,(R_1-R_D(p)+R_{\Omega}(p)-R_D(p))
$$
since $R_1\ge R_{\Omega}(p)$.  Thus theoretically (i) is better than (ii) since wo do not need to make use of 
the observation data beyond the time $2\sqrt{\rho/\mu}\,(R_1-R_D(p))$.

In the next subsection we give a proof of Theorem 4.1.

\subsection{Proof of Theorem 4.1}
We set $\mbox{\boldmath $R$}=\mbox{\boldmath $w$}-\mbox{\boldmath $w$}_0$
and $\Sigma=\Xi-\Xi_0$, where
$$\left\{
\begin{array}{ll}
\displaystyle
\Xi(x;\tau)=\int_0^Te^{-\tau t}\vartheta(x,t)\,dt, & x\in\Omega\setminus\overline D,\\
\\
\displaystyle
\Xi_0(x;\tau)=\int_0^Te^{-\tau t}\Theta(x,t)\,dt, & x\in\Bbb R^3.
\end{array}
\right.
$$ 
Note that $\Xi_0(x;\tau)\equiv 0$ since $\Theta(x,t)\equiv 0$.

It follows from (4.3) that $\nabla\cdot\mbox{\boldmath $w$}_0=0$. 
Applying this together with $\Theta=0$ to the equation (2.17) in \cite{IThermo}
and noting the bounds 
$$\displaystyle
\Vert\partial_t\mbox{\boldmath $u$}(\,\cdot,\,T)\Vert_{L^2(\Omega\setminus\overline D)}+
\Vert\mbox{\boldmath $u$}(\,\cdot,\,T)\Vert_{L^2(\Omega\setminus\overline D)}
+\Vert\nabla\cdot\mbox{\boldmath $u$}(\,\cdot,\,T)\Vert_{L^2(\Omega\setminus\overline D)}<\infty,
$$
$$\displaystyle
\Vert\partial_t\mbox{\boldmath $v$}_s(\,\cdot,\,T)\Vert_{L^2(\Omega)}
+\Vert\mbox{\boldmath $v$}_s(\,\cdot,\,T)\Vert_{L^2(\Omega)}<\infty,
$$
$$
\displaystyle
\Vert\vartheta(\,\cdot,\,T)\Vert_{L^2(\Omega\setminus\overline D)}<\infty,
$$
we have the expression
$$\begin{array}{ll}
\displaystyle
I^1(\tau;\mbox{\boldmath $v$}_s,0)
&
\displaystyle
=J(\tau)+
\left(E(\tau)+\frac{1}{\theta_0\tau}\,e(\tau)\right)\\
\\
\displaystyle
&
\,\,\,
\displaystyle
+O(e^{-\tau T}\left(\Vert \mbox{\boldmath $w$}_0\Vert_{L^2(\Omega)}+\Vert \mbox{\boldmath $R$}\Vert_{L^2(\Omega\setminus\overline D)}
+\Vert\Sigma\Vert_{L^2(\Omega\setminus\overline D)}\,\right)),
\end{array}
\tag {4.8}
$$
where
$$\left\{
\begin{array}
{l}
\displaystyle
J(\tau)=
\int_D\left(2\mu\vert\text{Sym}\,\nabla\mbox{\boldmath $w$}_0\vert^2+\rho\tau^2\vert\mbox{\boldmath $w$}_0\vert^2\right)\,dx,
\\
\\
\displaystyle
E(\tau)=\int_{\Omega\setminus\overline D}\left(2\mu\vert\text{Sym}\,\nabla\mbox{\boldmath $R$}\vert^2+\lambda\vert\nabla\cdot\mbox{\boldmath $R$}\vert^2
+\rho\tau^2\vert\mbox{\boldmath $R$}\vert^2\right)\,dx,
\\
\\
\displaystyle
e(\tau)=\int_{\Omega\setminus\overline D}\left(k\vert\nabla\Sigma\vert^2+c\tau\vert\Sigma\vert^2\right)\,dx.
\end{array}
\right.
$$

Since we have
$$\left\{
\begin{array}{l}
\displaystyle
\Vert\mbox{\boldmath $R$}\Vert_{L^2(\Omega\setminus\overline D)}\le\frac{1}{\sqrt{\rho}\tau}\,E(\tau)^{1/2},
\\
\\
\displaystyle
\Vert\Sigma\Vert_{L^2(\Omega\setminus\overline D)}\le\frac{\sqrt{\theta_0}}{\sqrt{c}}\cdot\frac{1}{\sqrt{\theta_0}\sqrt{\tau}}\,e(\tau)^{1/2},
\end{array}
\right.
$$
we obtain
$$\displaystyle
\Vert\mbox{\boldmath $R$}\Vert_{L^2(\Omega\setminus\overline D)}
+\Vert\Sigma\Vert_{L^2(\Omega\setminus\overline D)}
\le\frac{1}{2}\left(\frac{1}{\rho\tau^2}+\frac{\theta_0}{c}\right)+\frac{1}{2}\left(E(\tau)+\frac{1}{\theta_0\tau}e(\tau)\right).
$$
Thus, equation (4.8) becomes
$$\begin{array}{ll}
\displaystyle
I^1(\tau;\mbox{\boldmath $v$}_s,0)
&
\displaystyle
=J(\tau)+
(1+O(e^{-\tau T}))\left(E(\tau)+\frac{1}{\theta_0\tau}\,e(\tau)\right)+O(e^{-\tau T}\Vert \mbox{\boldmath $w$}_0\Vert_{L^2(\Omega)}).
\end{array}
\tag {4.9}
$$
By Proposition 2.4 in \cite{IThermo}, we have
$$\displaystyle
E(\tau)+\frac{1}{\theta_0\tau}e(\tau)=O(\tau^2 J(\tau)+\tau^2e^{-2\tau T}).
$$
Thus it follows from this and (4.9) that, for all $\tau\ge\tau_0>>1$,
$$\displaystyle
I^1(\tau;\mbox{\boldmath $v$}_s,0)\ge J(\tau)+O(e^{-\tau T}\Vert \mbox{\boldmath $w$}_0\Vert_{L^2(\Omega)})
\tag {4.10}
$$
and
$$\displaystyle
I^1(\tau;\mbox{\boldmath $v$}_s,0)=O(\tau^2J(\tau)+\tau^2e^{-2\tau T}+e^{-\tau T}\Vert \mbox{\boldmath $w$}_0\Vert_{L^2(\Omega)}).
\tag {4.11}
$$
Thus, everything is reduced to obtaining the upper and lower bounds for $J(\tau)$ and upper bound
for $\Vert\mbox{\boldmath $w$}_0\Vert_{L^2(\Omega)}$.

Let $\mbox{\boldmath $w$}_{s0}=\mbox{\boldmath $w$}_{s0}(\,\cdot\,;\tau)\in H^1(\Bbb R^3)$ solve
$$\begin{array}{ll}
\displaystyle
(\mu\Delta-\rho\tau^2)\mbox{\boldmath $w$}_{s0}+\rho\partial_t\mbox{\boldmath $v$}_{s}(x,0)=\mbox{\boldmath $0$},
& x\in\Bbb R^3.
\end{array}
\tag {4.12}
$$
We have the expression
$$\displaystyle
\mbox{\boldmath $w$}_{s0}(x;\tau)
=\frac{\rho}{\mu}\nabla\times\left\{
\frac{1}{4\pi}\int_{B_{R_2}(p)\setminus B_{R_1}(p)}\frac{e^{-\tau\sqrt{\rho/\mu}\,\vert x-y\vert}}{\vert x-y\vert}
\partial_t\mbox{\boldmath $\Phi$}(y,0)\,dy\,\right\}.
\tag {4.13}
$$
Since the function $\mbox{\boldmath $w$}_0$ satisfies the equation
$$\begin{array}
{ll}
\displaystyle
(\mu\Delta-\rho\tau^2)\mbox{\boldmath $w$}_0+\rho\partial_t\mbox{\boldmath $v$}_s(x,0)=\rho e^{-\tau T}(
\partial_t\mbox{\boldmath $v$}_s(x,T)+\tau\mbox{\boldmath $v$}_s(x,T)\,), & x\in\Bbb R^3,
\end{array}
$$
it follows from this together with (4.12) that we have
$$\displaystyle
\mbox{\boldmath $w$}_0(\,\cdot\,;\tau)=\mbox{\boldmath $w$}_{s0}(\,\cdot\,;\tau)+e^{-\tau T}\mbox{\boldmath $\epsilon$}_s,
$$
where the function $\mbox{\boldmath $\epsilon$}_s(\,\cdot,\,\tau)$ satisfies, as $\tau\rightarrow\infty$,
$$\displaystyle
\tau\Vert\mbox{\boldmath $\epsilon$}_s\Vert_{L^2(\Bbb R^3)}+\Vert\nabla\mbox{\boldmath $\epsilon$}_s\Vert_{L^2(\Bbb R^3)}
=O(1).
$$

Thus, one gets
$$
\left\{
\begin{array}{l}
\displaystyle
\Vert\mbox{\boldmath $w$}_0\Vert_{L^2(\Omega)}\le \Vert \mbox{\boldmath $w$}_{s0}\Vert_{L^2(\Omega)}+O(\tau^{-1}e^{-\tau T}),
\\
\\
\displaystyle
J(\tau)\ge \rho\tau^2\Vert \mbox{\boldmath $w$}_{s0}\Vert_{L^2(D)}^2+O(e^{-2\tau T}),\\
\\
\displaystyle
J(\tau)=O(\Vert\nabla\mbox{\boldmath $w$}_{s0}\Vert_{L^2(D)}^2+\tau^2\Vert\mbox{\boldmath $w$}_{s0}\Vert_{L^2(D)}^2+e^{-2\tau T}).
\end{array}
\right.
$$
Applying these to (4.10) and (4.11), we obtain
$$\displaystyle
I^1(\tau;\mbox{\boldmath $v$}_s,0)\ge \rho\tau^2\Vert \mbox{\boldmath $w$}_{s0}\Vert_{L^2(D)}^2
+O(e^{-\tau T}\Vert \mbox{\boldmath $w$}_{s0}\Vert_{L^2(\Omega)})+O(e^{-2\tau T})
\tag {4.14}
$$
and
$$\displaystyle
I^1(\tau;\mbox{\boldmath $v$}_s,0)=O(\tau^2\Vert\nabla\mbox{\boldmath $w$}_{s0}\Vert_{L^2(D)}^2
+\tau^4\Vert\mbox{\boldmath $w$}_{s0}\Vert_{L^2(D)}^2
+\tau^2e^{-2\tau T}+e^{-\tau T}\Vert \mbox{\boldmath $w$}_{s0}\Vert_{L^2(\Omega)}).
\tag {4.15}
$$

\proclaim{\noindent Lemma 4.2.}
We have, for all $x\in B_{R_1}(p)$,
$$
\displaystyle
\mbox{\boldmath $w$}_{s0}(x;\tau)=
\frac{\rho}{\mu}M(\tau\sqrt{\rho/\mu};R_1,R_2)
e^{-\tau\sqrt{\rho/\mu}\,R_1}
\nabla\left(\frac{\sinh\tau\sqrt{\rho/\mu}\vert x-p\vert}{\vert x-p\vert}\right)\times\mbox{\boldmath $a$},
$$
where the coefficient $M(\tau\sqrt{\rho/\mu};R_1,R_2)$ is independent of $x$ and satifies
$$\displaystyle
\lim_{\tau\rightarrow\infty}\tau^4M(\tau;R_1,R_2)=2R_1(R_1-R_2)^2.
$$
\endproclaim

\noindent
Since $\Omega\subset B_{R_{\Omega}(p)}(p)$, from Lemma 4.2 we have immediately the following estimate:
$$\displaystyle
\Vert \mbox{\boldmath $w$}_{s0}\Vert_{L^2(\Omega)}=O(\tau^{-3} e^{-\tau\sqrt{\rho/\mu}(R_1-R_{D}(p))}
e^{\tau\sqrt{\rho/\mu}(R_{\Omega}(p)-R_D(p))})
\tag {4.16}
$$
and thus
$$\displaystyle
e^{-\tau T}\Vert \mbox{\boldmath $w$}_{s0}\Vert_{L^2(\Omega)}=O(\tau^{-3} e^{-\tau\sqrt{\rho/\mu}(R_1-R_{D}(p))}
e^{-\tau\{T-\sqrt{\rho/\mu}\,(R_{\Omega}(p)-R_D(p))\}})
$$
or
$$\begin{array}{ll}
\displaystyle
e^{2\tau\sqrt{\rho/\mu}(R_1-R_D(p))}e^{-\tau T}\Vert \mbox{\boldmath $w$}_{s0}\Vert_{L^2(\Omega)}
&
\displaystyle
=O(\tau^{-3} e^{\tau\sqrt{\rho/\mu}(R_1-R_{D}(p))}
e^{-\tau\{T-\sqrt{\rho/\mu}\,(R_{\Omega}(p)-R_D(p))\}})
\\
\\
\displaystyle
&
\displaystyle
=O(\tau^{-3}
e^{-\tau\{T-\sqrt{\rho/\mu}\,(R_1-R_D(p)+R_{\Omega}(p)-R_D(p))\}}).
\end{array}
\tag {4.17}
$$

On $D$ we have also
$$\displaystyle
\Vert\nabla\mbox{\boldmath $w$}_{s0}\Vert_{L^2(D)}^2+\tau^2\Vert\mbox{\boldmath $w$}_{s0}\Vert_{L^2(D)}^2
=O(\tau^{-4} e^{-2\tau\sqrt{\rho/\mu}(R_1-R_{D}(p))}).
\tag {4.18}
$$
On the other hand, the lower estimate of $\Vert\mbox{\boldmath $w$}_{s0}\Vert_{L^2(D)}$ is
not so trivial.

\proclaim{\noindent Lemma 4.3.}
Assume that $\partial D$ is $C^2$.  Then there exist positive numbers $\tau_0$, $C$ and nonnegative number $\alpha$ such that,
for all $\tau\ge\tau_0$
$$\displaystyle
\tau^{2\alpha}e^{2\tau\sqrt{\rho/\mu}\,(R_1-R_D(p))}\,\Vert\mbox{\boldmath $w$}_{s0}\Vert_{L^2(D)}^2\ge C.
$$

\endproclaim

\noindent
From (4.14), (4.17) and Lemma 4.3 we obtain
$$\begin{array}{ll}
\displaystyle
e^{2\tau\sqrt{\rho/\mu}(R_1-R_D(p))}I^1(\tau;\mbox{\boldmath $v$}_s,0)
&
\displaystyle
\ge \rho\tau^{2-2\alpha}\cdot\tau^{2\alpha}e^{2\tau\sqrt{\rho/\mu}(R_1-R_D(p))}\Vert \mbox{\boldmath $w$}_{s0}\Vert_{L^2(D)}^2
\\
\\
\displaystyle
&
\displaystyle
\,\,\,
+O(e^{2\tau\sqrt{\rho/\mu}(R_1-R_D(p))}\,e^{-\tau T}\Vert \mbox{\boldmath $w$}_{s0}\Vert_{L^2(\Omega)})
\\
\\
\displaystyle
&
\displaystyle
\,\,\,
+O(e^{-\tau T}e^{-\tau\{T-2\sqrt{\rho/\mu}(R_1-R_D(p))\}})
\\
\\
\displaystyle
&
\displaystyle
\ge \rho\tau^{2-2\alpha}C
\\
\\
\displaystyle
&
\displaystyle
\,\,\,
+O(\tau^{-3}
e^{-\tau\{T-\sqrt{\rho/\mu}\,(R_1-R_D(p)+R_{\Omega}(p)-R_D(p))\}})
\\
\\
\displaystyle
&
\displaystyle
\displaystyle
\,\,\,
+O(e^{-\tau T}e^{-\tau\{T-2\sqrt{\rho/\mu}(R_1-R_D(p))\}}).
\end{array}
$$
Thus, if $T$ satisfies (4.5), then $T>2\sqrt{\rho/\mu}(R_1-R_D(p))$ and one gets
$$\displaystyle
\liminf_{\tau\rightarrow\infty}\tau^{2\alpha-2}e^{2\tau\sqrt{\rho/\mu}(R_1-R_D(p))}I^1(\tau;\mbox{\boldmath $v$}_s,0)>0.
\tag {4.19}
$$
Besides, from (4.15), (4.17) and (4.18) we obtain
$$\begin{array}{ll}
\displaystyle
e^{2\tau\sqrt{\rho/\mu}(R_1-R_D(p))}I^1(\tau;\mbox{\boldmath $v$}_s,0)
&
\displaystyle
=O(\tau^{-2}+\tau^2 e^{-2\tau (T-\sqrt{\rho/\mu}(R_1-R_D(p))})
\\
\\
\displaystyle
&
\displaystyle
\,\,\,
+O(\tau^{-3}
e^{-\tau\{T-\sqrt{\rho/\mu}\,(R_1-R_D(p)+R_{\Omega}(p)-R_D(p))\}}).
\end{array}
$$
Thus if $T$ satisfies
$$\displaystyle
T\ge \sqrt{\frac{\rho}{\mu}}\,(R_1-R_D(p)+R_{\Omega}(p)-R_D(p)),
$$
then we obtain
$$\displaystyle
\limsup_{\tau\rightarrow\infty}\tau^2e^{2\tau\sqrt{\rho/\mu}(R_1-R_D(p))}I^1(\tau;\mbox{\boldmath $v$}_s,0)<\infty.
\tag {4.20}
$$
Now Theorem 4.1 (i) is a direct consequence of (4.19) and (4.20).

Next we give a proof of (4.7).
From (4.14), (4.16) and Lemma 4.3 we have
$$\begin{array}{ll}
\displaystyle
e^{\tau T}I^1(\tau;\mbox{\boldmath $v$}_s,0)
&
\ge \rho\,C\tau^{2-2\alpha}e^{\tau\{T-2\sqrt{\rho/\mu}(R_1-R_D(p))\}}
\\
\\
\displaystyle
&
\,\,\,
\displaystyle
+O(\tau^{-3} e^{-\tau\sqrt{\rho/\mu}(R_1-R_{\Omega}(p))})+O(e^{-\tau T}).
\end{array}
$$
Thus if $T>2\sqrt{\rho/\mu}(R_1-R_D(p))$, then we have 
$$\displaystyle
\lim_{\tau\rightarrow\infty}e^{\tau T}I^1(\tau;\mbox{\boldmath $v$}_s,0)=\infty.
$$
From (4.15), (4.16) and (4.18), we obtain
$$\begin{array}{ll}
\displaystyle
e^{\tau T}I^1(\tau;\mbox{\boldmath $v$}_s,0)
&
\displaystyle
=O(\tau^{-2}e^{\tau\{T-2\sqrt{\rho/\mu}(R_1-R_D(p))\}}e^{-\tau\sqrt{\rho/\mu}(R_1-R_D(p))})
+O(\tau^2e^{-\tau T})
\\
\\
\displaystyle
&
\displaystyle
\,\,\,
+O(\tau^{-3} e^{-\tau\sqrt{\rho/\mu}(R_1-R_{\Omega}(p))}).
\end{array}
$$
Thus, if $T\le 2\sqrt{\rho/\mu}(R_1-R_D(p))$, then we have 
$$\displaystyle
\lim_{\tau\rightarrow\infty}e^{\tau T}I^1(\tau;\mbox{\boldmath $v$}_s,0)=0.
$$
This completes the proof of Theorem 4.1.

The proofs of Lemmas 4.2 and 4.3 are given in Sections 5 and 6, respectively.

\section{Proof of Lemma 4.2}

We derive a recurrence formula for the sequence of the integrals
$$\begin{array}{ll}
\displaystyle
v_j(x)=\int_B\frac{e^{-\tau\vert x-y\vert}}{\vert x-y\vert}\,\vert y\vert^j\,dy, & x\in B,
\end{array}
$$
where $B=\{y\in\Bbb R^3\,\vert\,\vert y\vert<\eta\}$, with $\eta>0$ and $j=-1,0,1,2,\cdots$.

\proclaim{\noindent Proposition 5.1.}  Let $j\ge 1$.  We have
$$
\displaystyle
v_j(x)
\displaystyle
=\frac{j(j+1)}{\tau^2}v_{j-2}(x)
+
\frac{4\pi}{\vert x\vert\,\tau^2}
\left\{
\vert x\vert^{j+1}-\left(\eta+\frac{j+1}{\tau}\right)\eta^je^{-\tau\eta}\sinh\tau\vert x\vert
\right\}.
\tag {5.1}
$$

\endproclaim

{\it\noindent Proof.}
By (A.1) in \cite{EV}, we have
$$\displaystyle
v_j(x)=\frac{2\pi}{\vert x\vert\,\tau}K_j(\vert x\vert),
$$
where
$$\displaystyle
K_j(\xi)=\int_0^{\eta}(e^{-\tau\vert\xi-r\vert}-e^{-\tau(\xi+r)})\,r^{j+1}\,dr.
$$
Thus, it suffices to prove that
$$\displaystyle
K_j(\xi)=\frac{j(j+1)}{\tau^2}K_{j-2}(\xi)+\frac{2}{\tau}
\left\{\xi^{j+1}-\left(\eta+\frac{j+1}{\tau}\right)\eta^je^{-\tau\eta}\sinh\tau\xi
\right\}.
\tag {5.2}
$$
First write
$$\displaystyle
K_j(\xi)=2e^{-\tau\xi}\int_0^{\xi}r^{j+1}\sinh\tau r\,dr
+2\sinh\tau\xi\int_{\xi}^{\eta}r^{j+1}e^{-\tau r}\,dr.
$$
Integration by parts yields
$$\begin{array}{l}
\displaystyle
\,\,\,\,\,\,K_j(\xi)
\\
\\
\displaystyle
=
\frac{2e^{-\tau\xi}}{\tau}\left\{\xi^{j+1}\cosh\tau\xi-(j+1)\int_0^{\xi}r^j\cosh\tau r\,dr\right\}
\\
\\
\displaystyle
\,\,\,
+\frac{2\sinh\tau\xi}{\tau}
\left\{\xi^{j+1}e^{-\tau\xi}-\eta^{j+1}e^{-\tau\eta}+(j+1)\int_{\xi}^{\eta}r^je^{-\tau r}\,dr\right\}
\\
\\
\displaystyle
=\frac{2}{\tau}(\xi^{j+1}-\eta^{j+1}e^{-\tau\eta}\sinh\tau\xi)+\frac{2(j+1)}{\tau}L_j(\xi),
\end{array}
\tag {5.3}
$$
where
$$\displaystyle
L_j(\xi)=-e^{-\tau\xi}\int_0^{\xi}r^j\cosh\tau r\,dr+\sinh\tau\xi\int_{\xi}^{\eta}r^je^{-\tau r}\,dr.
$$
Integration by parts yields
$$\begin{array}{l}
\displaystyle
\,\,\,\,\,\,
L_j(\xi)
\\
\\
\displaystyle
=-\frac{e^{-\tau\xi}}{\tau}
\left(\xi^j\sinh\tau\xi-j\int_0^{\xi}r^{j-1}\sinh\tau r\,dr\right)
\\
\\
\displaystyle
\,\,\,
-\frac{\sinh\tau\xi}{\tau}
\left(\eta^je^{-\tau\eta}-\xi^je^{-\tau\xi}-j\int_{\xi}^{\eta}r^{j-1}e^{-\tau r}\,dr\right)\\
\\
\displaystyle
=-\frac{\eta^j}{\tau}e^{-\tau\eta}\sinh\tau\xi+\frac{j}{2\tau}K_{j-2}.
\end{array}
$$
Substituting this into (5.3), we obtain (5.2).

\noindent
$\Box$

\noindent
Next we compute
$$\begin{array}{ll}
\displaystyle
J(x)=\frac{1}{4\pi}\int_{B_2\setminus B_1}\frac{e^{-\tau\vert x-y\vert}}{\vert x-y\vert}\,
(R_1-\vert y-p\vert)^2(R_2-\vert y-p\vert)^2\,dy, & x\in B_1,
\end{array}
$$
where $B_j=B_{R_j}(p)$.  We have 
$$\displaystyle
J(x)=J_2(x)-J_1(x),
$$
where
$$\begin{array}{ll}
\displaystyle
J_i(x)=\frac{1}{4\pi}\int_{B_i}\frac{e^{-\tau\vert x-y\vert}}{\vert x-y\vert}\,
(R_1-\vert y-p\vert)^2(R_2-\vert y-p\vert)^2\,dy, & i=1,2.
\end{array}
$$
Set
$$\begin{array}{lll}
\displaystyle
v_j^i(x)=\frac{1}{4\pi}\int_{B_i}\frac{e^{-\tau\vert x-y\vert}}{\vert x-y\vert}\,\vert y-p\vert^j\,dy, & j\ge1, & i=1,2.
\end{array}
$$
One has the expression
$$\begin{array}{ll}
\displaystyle
J_i(x)
&
\displaystyle
=(R_1R_2)^2v_0^i(x)-2R_1R_2(R_1+R_2)v_1^i
\\
\\
\displaystyle
&
\displaystyle
\,\,\,
+\left\{(R_1+R_2)^2+2R_1R_2\right\}v_2^i
-2(R_1+R_2)v_3^i+v_4^i.
\end{array}
$$
Thus one gets
$$\begin{array}{ll}
\displaystyle
J(x)
&
\displaystyle
=(R_1R_2)^2(v_0^2(x)-v_0^1(x))-2R_1R_2(R_1+R_2)(v_1^2(x)-v_1^1(x))
\\
\\
\displaystyle
&
\displaystyle
\,\,\,
+\left\{(R_1+R_2)^2+2R_1R_2\right\}(v_2^2(x)-v_2^1(x))
\\
\\
\displaystyle
&
\displaystyle
\,\,\,
-2(R_1+R_2)(v_3^2(x)-v_3^1(x))+(v_4^2(x)-v_4^1(x)).
\end{array}
$$
Here from (5.1) we have
$$\begin{array}{ll}
\displaystyle
\tau^2(v_3^2(x)-v_3^1(x))
&
\displaystyle
=12(v_1^2(x)-v_1^1(x))
\\
\\
\displaystyle
&
\displaystyle
\,\,\,
-\left\{\left(R_2+\frac{4}{\tau}\right)R_2^3e^{-\tau R_2}
-\left(R_1+\frac{4}{\tau}\right)R_1^3e^{-\tau R_1}\right\}
\frac{\sinh\tau\vert x\vert}{\vert x\vert}
\\
\\
\displaystyle
&
\displaystyle
=12(v_1^2(x)-v_1^1(x))
\\
\\
\displaystyle
&
\displaystyle
\,\,\,
+\left\{\left(R_1+\frac{4}{\tau}\right)R_1^3
+O(e^{-\tau(R_2-R_1)})\right\}e^{-\tau R_1}
\frac{\sinh\tau\vert x\vert}{\vert x\vert}
\end{array}
$$
and
$$\begin{array}{ll}
\displaystyle
\tau^2(v_4^2(x)-v_4^1(x))
&
\displaystyle
=20(v_2^2(x)-v_2^1(x))
\\
\\
\displaystyle
&
\displaystyle
\,\,\,
-\left\{\left(R_2+\frac{5}{\tau}\right)R_2^4e^{-\tau R_2}
-\left(R_1+\frac{5}{\tau}\right)R_1^4e^{-\tau R_1}\right\}
\frac{\sinh\tau\vert x\vert}{\vert x\vert}
\\
\\
\displaystyle
&
\displaystyle
=20(v_2^2(x)-v_2^1(x))
\\
\\
\displaystyle
&
\displaystyle
\,\,\,
+\left\{\left(R_1+\frac{5}{\tau}\right)R_1^4+O(e^{-\tau(R_2-R_1)})\right\}e^{-\tau R_1}
\frac{\sinh\tau\vert x\vert}{\vert x\vert}.
\end{array}
$$

In the proof of Lemma 2.2 we have already shown that
$$
\left\{
\begin{array}{l}
\displaystyle
\tau^2(v_0^2(x)-v_0^1(x))
=\left(R_1+\frac{1}{\tau}+O(e^{-\tau (R_2-R_1)})\right)e^{-\tau R_1}
\frac{\sinh\tau\vert x-p\vert}{\vert x-p\vert},
\\
\\
\displaystyle
\tau^2(v_1^2(x)-v_1^1(x))
=\left(R_1^2+\frac{2}{\tau}R_1+\frac{2}{\tau^2}
+O(e^{-\tau(R_2-R_1)})
\right)
e^{-\tau R_1}\frac{\sinh\tau\vert x-p\vert}{\vert x-p\vert},
\\
\\
\displaystyle
\tau^2(v_2^2(x)-v_2^1(x))
=\left(
R_1^3+\frac{3R_1^2}{\tau}+\frac{6R_1}{\tau^2}+\frac{6}{\tau^3}
+O(e^{-\tau(R_2-R_1)})
\right)
e^{-\tau R_1}
\frac{\sinh\tau\vert x-p\vert}{\vert x-p\vert}.
\end{array}
\right.
$$
From these we obtain
$$\begin{array}{ll}
\displaystyle
J(x)
&
\displaystyle
=(R_1R_2)^2(v_0^2(x)-v_0^1(x))-2R_1R_2(R_1+R_2)(v_1^2(x)-v_1^1(x))
\\
\\
\displaystyle
&
\displaystyle
\,\,\,
+\left\{(R_1+R_2)^2+2R_1R_2\right\}(v_2^2(x)-v_2^1(x))
\\
\\
\displaystyle
&
\displaystyle
\,\,\,
-2(R_1+R_2)(v_3^2(x)-v_3^1(x))+(v_4^2(x)-v_4^1(x))
\\
\\
\displaystyle
&
\displaystyle
=
M(\tau;R_1,R_2)e^{-\tau R_1}
\frac{\sinh\tau\vert x-p\vert}{\vert x-p\vert},
\end{array}
$$
where
$$\begin{array}{ll}
\displaystyle
M(\tau;R_1,R_2)
&
\displaystyle
=A(R_1,R_2)\tau^{-2}+B(R_1,R_2)\tau^{-3}+C(R_1,R_2)\tau^{-4}+O(\tau^{-5}),
\end{array}
$$
$$\begin{array}{ll}
\displaystyle
\,\,\,\,\,\,
A(R_1,R_2)
\\
\\
\displaystyle
=R_1^3\left\{R_2^2-2R_2(R_1+R_2)+(R_1+R_2)^2+2R_1R_2-2(R_1+R_2)R_1+R_1^2\right\}
\\
\\
\displaystyle
=R_1^3\{R_2^2-2R_2^2+(R_1+R_2)^2-2(R_1+R_2)R_1+R_1^2\}
\\
\\
\displaystyle
=R_1^3(-R_2^2+R_2^2)
\\
\\
\displaystyle
=0,
\end{array}
$$
$$\begin{array}{l}
\displaystyle
\,\,\,\,\,\,
B(R_1,R_2)
\\
\\
\displaystyle
=R_1^2\left\{R_2^2-4R_2(R_1+R_2)+3(R_1+R_2)^2+6R_1R_2-8(R_1+R_2)R_1+5R_1^2\right\}
\\
\\
\displaystyle
=R_1^2\left\{R_2^2-4R_2^2+2R_1R_2+3(R_1+R_2)^2-8R_2R_1-3R_1^2\right\}
\\
\\
\displaystyle
=R_1^2(-3R_2^2-6R_1R_2+6R_1R_2+3R_2^2)
\\
\\
\displaystyle
=0
\end{array}
$$
and
$$\begin{array}{l}
\displaystyle
\,\,\,\,\,\,
C(R_1,R_2)
\\
\\
\displaystyle
=-4R_1R_2(R_1+R_2)+6R_1\{(R_1+R_2)^2+2R_1R_2\}-24(R_1+R_2)R_1^2+20R_1^3
\\
\\
\displaystyle
=(R_1+R_2)\left\{6R_1(R_1+R_2)-4R_1R_2-24R_1^2\right\}+12R_1^2R_2+20R_1^3
\\
\\
\displaystyle
=R_1(R_1+R_2)\left\{6(R_1+R_2)-4R_2-24R_1\right\}+4R_1^2(3R_2+5R_1)
\\
\\
\displaystyle
=2R_1(R_1+R_2)(R_2-9R_1)+4R_1^2(3R_2+5R_1)\\
\\
\displaystyle
=2R_1\left\{(R_1+R_2)(R_2-9R_1)+2R_1(3R_2+5R_1)\right\}
\\
\\
\displaystyle
=2R_1(R_1R_2-9R_1^2+R_2^2-9R_1R_2+6R_1R_2+10R_1^2)
\\
\\
\displaystyle
=2R_1(R_1-R_2)^2.
\end{array}
$$

Thus, one gets
$$\displaystyle
\lim_{\tau\rightarrow\infty}\tau^4M(\tau;R_1,R_2)=2R_1(R_1-R_2)^2.
$$

From these together with the third equation on (4.2) and (4.13) we obtain the desired conclusion.

\section{Proof of Lemma 4.3}

Choose a point $q\in\partial D$ such that $\vert q-p\vert=R_D(p)$.  We have $\nu(q)=\frac{q-p}{\vert q-p\vert}$.
Since $\partial D$ is $C^2$, one can find an open ball $B=B_{\delta}(q-\delta \nu(q))$ with a small $\delta>0$
such that $B\subset D$, $\partial B\cap\partial D=\{q\}$, $R_B(p)=R_D(p)$ and $p\in\Bbb R^3\setminus\overline B$.
Note that the center point $q-\delta\nu(q)$ is on the segment $[p,q]$.

Since we have $\Vert \mbox{\boldmath $w$}_{s0}\Vert_{L^2(D)}\ge\Vert \mbox{\boldmath $w$}_{s0}\Vert_{L^2(B)}$, it suffices to 
prove Lemma 4.3 in the case when $D$ is replaced with the $B$.
To this end, let $0<\delta'<\delta$ and $B'=B_{R_D(p)-\delta'}(p)$.  We parametrize the set $B\setminus\overline B'$ by using the polar
coordinates centered at $p$.  Let $\mbox{\boldmath $\omega$}\in S^2$ and $R_D(p)-\delta'<r<R_D(p)$.

The point
$x(r,\mbox{\boldmath $\omega$})\equiv p+r\mbox{\boldmath $\omega$}$ 
belongs to set $B\setminus\overline B'$
if and only if $x(r,\mbox{\boldmath $\omega$})\in B$, that is 
$$\displaystyle
\vert x(r,\mbox{\boldmath $\omega$})-(q-\delta\nu(q))\vert<\delta.
$$
Since $q=p+R_D(p)\nu(q)$, this is equivalent to the inequality
%$$\displaystyle
%\vert r\mbox{\boldmath $\omega$}-(R_D(p)-\delta)\nu(q)\vert<\delta.
%$$
%$$\displaystyle
%r^2-2r(R_D(p)-\delta)\mbox{\boldmath $\omega$}\cdot\nu(q)+(R_D(p)-\delta)^2<\delta^2.
%$$
$$\displaystyle
\mbox{\boldmath $\omega$}\cdot\nu(q)
>\frac{r^2+(R_D(p)-\delta)^2-\delta^2}{2r(R_D(p)-\delta)}.
$$
Thus, we have the expression
$$\displaystyle
B\setminus\overline B'=\cup_{R_D(p)-\delta'<r<R_D(p)}\,S(r),
$$
where
$$\displaystyle
S(r)=\left\{\mbox{\boldmath $\omega$}\in S^2\,\vert\,\mbox{\boldmath $\omega$}\cdot\nu(q)
>\frac{r^2+(R_D(p)-\delta)^2-\delta^2}{2r(R_D(p)-\delta)}\right\}.
$$

Next we parametrize the surface $S(r)$ in the following way.

Let $\theta=\theta(r)\in \,]0,\,\frac{\pi}{2}[$ denote the unique solution of the equation
$$\displaystyle
\cos\theta=\frac{r^2+(R_D(p)-\delta)^2-\delta^2}{2r(R_D(p)-\delta)}.
\tag {6.1}
$$
Choose two linearly independent vectors $\mbox{\boldmath $b$}$ and $\mbox{\boldmath $c$}$ in such a way that
$\mbox{\boldmath $b$}\cdot\mbox{\boldmath $c$}=0$ and $\mbox{\boldmath $b$}\times\mbox{\boldmath $c$}=\nu(q)$.
Let $R_{D}(p)-\delta'<r<R_D(p)$, $0<s<r\sin\theta(r)$, $0\le\phi<2\pi$ and $0\le h$.
Set
$$\displaystyle
\Upsilon(r,s,\phi)=p+r\nu(q)+s(\cos\phi\mbox{\boldmath $b$}+\sin\phi\mbox{\boldmath $c$})-h\nu(q).
$$
Here, the $h$ is {\it unknown} and determined by the equation
$$\displaystyle
\vert
\Upsilon(r,s,\phi)-p\vert=r.
$$
Solving this equation and choosing the smaller one, we have
$$\displaystyle
h=r-\sqrt{r^2-s^2}.
$$
Thus we have
$$\displaystyle
\Upsilon(r,s,\phi)=p+\sqrt{r^2-s^2}\nu(q)+s(\cos\phi\mbox{\boldmath $b$}+\sin\phi\mbox{\boldmath $c$})
$$
and using this, we have the expression
$$\displaystyle
S(r)\setminus\{p+r\nu(q)\}=\left\{\Upsilon(r,s,\phi)\,\vert\,0<s<r\sin\theta(r), 0\le\phi<2\pi\right\}.
$$

Define
$$\displaystyle
G=\left\{(r,s,\phi)\,\vert\,R_D(p)-\delta'<r<R_D(p),\,0<s<r\sin\theta(r),\,0<\phi<2\pi\right\}.
$$
We see that the map
$$\displaystyle
\Upsilon:G\ni(r,s,\phi)\mapsto\Upsilon(r,s,\phi)\in (B\setminus\overline B')\setminus Z,
$$
is bijective, where 
$$
\displaystyle
Z=\left\{p+\sqrt{r^2-s^2}\nu(q)+s\mbox{\boldmath $b$}\,\vert, R_D(p)-\delta'<r<R_D(p),\,0\le s<r\sin\theta(r)\right\}.
$$
Note that the Lebesgue measure of $Z$ is zero.

We have
$$\displaystyle
\Upsilon'(r,s,\phi)
=\left(
\begin{array}{ccc}
\displaystyle
\frac{r}{\sqrt{r^2-s^2}}\nu(q) & 
\displaystyle
-\frac{s}{\sqrt{r^2-s^2}}\nu(q)+\cos\phi\mbox{\boldmath $b$}+\sin\phi\mbox{\boldmath $c$}
&
\displaystyle
s(-\sin\phi\mbox{\boldmath $b$}+\cos\phi\mbox{\boldmath $c$})
\end{array}
\right)
$$
and thus
$$\displaystyle
\text{det}\,\Upsilon'(r,s,\phi)=\frac{sr}{\sqrt{r^2-s^2}}.
$$
Here we have
$$\displaystyle
\nabla
\left(\frac{\sinh\tau\vert x-p\vert}{\vert x-p\vert}\,\right)\times\mbox{\boldmath $a$}
=\left(\tau\frac{\cosh\tau\vert x-p\vert}{\vert x-p\vert}-\frac{\sinh\tau\vert x-p\vert}{\vert x-p\vert^2}\right)\,
\frac{x-p}{\vert x-p\vert}\times\mbox{\boldmath $a$}.
$$
Let $x\in B\setminus\overline B'$.  Since we have $R_D(p)-\delta'<\vert x-p\vert<R_D(p)$, one gets
$$\displaystyle
\left\vert\tau\frac{\cosh\tau\vert x-p\vert}{\vert x-p\vert}-\frac{\sinh\tau\vert x-p\vert}{\vert x-p\vert^2}\right\vert
\ge C\tau e^{\tau\vert x-p\vert},
$$
where $C$ is a positive number being independent of $\tau>>1$ and $x$.
Thus one has
$$\displaystyle
\int_B\left\vert\nabla\left(\frac{\sinh\tau\vert x-p\vert}{\vert x-p\vert}\,\right)\times\mbox{\boldmath $a$}\right\vert^2\,dx
\ge C^2\tau^2\int_{(B\setminus\overline B')\setminus Z}e^{2\tau\vert x-p\vert}
\left\vert\frac{x-p}{\vert x-p\vert}\times\mbox{\boldmath $a$}\right\vert^2\,dx.
\tag {6.2}
$$

Now we compute the integral
$$\displaystyle
I(\tau)\equiv \int_{(B\setminus\overline B')\setminus Z}e^{2\tau\vert x-p\vert}
\left\vert\frac{x-p}{\vert x-p\vert}\times\mbox{\boldmath $a$}\right\vert^2\,dx.
$$
The change of variables $x=\Upsilon(r,s,\phi)$ yields
$$\displaystyle
I(\tau)=\int_{R_D(p)-\delta'}^{R_D(p)}dr\int_0^{r\sin\theta(r)}ds\int_0^{2\pi}d\phi
\frac{sr}{\sqrt{r^2-s^2}}e^{2\tau r}r^{-2}\left\vert(\Upsilon(r,s,\phi)-p)\times\mbox{\boldmath $a$}\right\vert^2.
$$
Since we have
$$\displaystyle
(\Upsilon(r,s,\phi)-p)\times\mbox{\boldmath $a$}
=\sqrt{r^2-s^2}\nu(q)\times\mbox{\boldmath $a$}
+s\cos\phi\,\mbox{\boldmath $b$}\times\mbox{\boldmath $a$}
+s\sin\phi\,\mbox{\boldmath $c$}\times\mbox{\boldmath $a$},
$$
one gets
$$\begin{array}{l}
\displaystyle
\,\,\,\,\,\,
\int_0^{2\pi}\vert(\Upsilon(r,s,\phi)-p)\times\mbox{\boldmath $a$}\vert^2\,d\phi
\\
\\
\displaystyle
=2\pi (r^2-s^2)\vert\nu(q)\times\mbox{\boldmath $a$}\vert^2
+s^2\int_0^{2\pi}\cos^2\phi\,d\phi\vert\mbox{\boldmath $b$}\times\mbox{\boldmath $a$}\vert^2
+s^2\int_0^{2\pi}\sin^2\phi\,d\phi\vert\mbox{\boldmath $c$}\times\mbox{\boldmath $a$}\vert^2
\\
\\
\displaystyle
=2\pi (r^2-s^2)\vert\nu(q)\times\mbox{\boldmath $a$}\vert^2
+\pi s^2(\vert\mbox{\boldmath $b$}\times\mbox{\boldmath $a$}\vert^2+
\vert\mbox{\boldmath $c$}\times\mbox{\boldmath $a$}\vert^2).
\end{array}
$$
Thus we obtain
$$\begin{array}{ll}
\displaystyle
I(\tau)
&
\displaystyle
=2\pi\int_{R_D(p)-\delta'}^{R_D(p)}dr\int_0^{r\sin\theta(r)}ds
\frac{s\sqrt{r^2-s^2}}{r}e^{2\tau r}\,\vert\nu(q)\times\mbox{\boldmath $a$}\vert^2
\\
\\
\displaystyle
&
\displaystyle
\,\,\,
+\pi\int_{R_D(p)-\delta'}^{R_D(p)}dr\int_0^{r\sin\theta(r)}ds
\frac{s^3}{r\sqrt{r^2-s^2}}\,
e^{2\tau r}(\vert\mbox{\boldmath $b$}\times\mbox{\boldmath $a$}\vert^2+
\vert\mbox{\boldmath $c$}\times\mbox{\boldmath $a$}\vert^2).
\end{array}
$$
Here we have
$$\left\{
\begin{array}{l}
\displaystyle
\int_0^{r\sin\theta(r)}s\sqrt{r^2-s^2}\,ds
=\frac{r^3}{3}(1-\cos^3\theta(r)),
\\
\\
\displaystyle
\int_0^{r\sin\theta(r)}\frac{s^3}{\sqrt{r^2-s^2}}\,ds
=r^3
\left\{(1-\cos\theta(r))-\frac{1}{3}(1-\cos^3\theta(r))\right\}.
\end{array}
\right.
$$
From (6.1) we have
$$\begin{array}{ll}
\displaystyle
1-\cos\theta(r)
&
\displaystyle
=\frac{(R_D(p)-r)(r-R_D(p)+2\delta)}{2r(R_D(p)-\delta)}
\\
\\
\displaystyle
&
\displaystyle
\ge (R_D(p)-r)\frac{2(\delta-\delta')}{2R_D(p)(R_D(p)-\delta)}
\\
\\
\displaystyle
&
\displaystyle
\equiv C'(R_D(p)-r)
\end{array}
$$
and
$$\begin{array}{l}
\displaystyle
\,\,\,\,\,\,
(1-\cos\theta(r))-\frac{1}{3}(1-\cos^3\theta(r))
\\
\\
\displaystyle
=(1-\cos\theta(r))\left\{1-\frac{1}{3}(1+\cos\theta(r)+\cos^2\theta(r))\right\}
\\
\\
\displaystyle
=\frac{1}{3}(1-\cos\theta(r))^2(2+\cos\theta(r))
\\
\\
\displaystyle
\ge \frac{2(C')^2}{3}(R_D(p)-r)^2.
\end{array}
$$
From these we obtain
$$\begin{array}{ll}
\displaystyle
I(\tau)
&
\displaystyle
\ge C_1
\int_{R_D(p)-\delta'}^{R_D(p)}(R_D(p)-r)e^{2\tau r}\,dr\vert\nu(q)\times\mbox{\boldmath $a$}\vert^2
\\
\\
\displaystyle
&
\displaystyle
\,\,\,
+C_2\int_{R_D(p)-\delta'}^{R_D(p)}(R_D(p)-r)^2e^{2\tau r}\,dr(\vert\mbox{\boldmath $b$}\times\mbox{\boldmath $a$}\vert^2+
\vert\mbox{\boldmath $c$}\times\mbox{\boldmath $a$}\vert^2).
\end{array}
$$
We see that, as $\tau\rightarrow\infty$
$$\displaystyle
\int_{R_D(p)-\delta'}^{R_D(p)}(R_D(p)-r)^j\,e^{2\tau r}\,dr
\sim
\frac{e^{2\tau R_D(p)}}{\tau^{j+1}}\int_0^{\infty}r^je^{-2r}\,dr.
$$
Thus, we conclude that: 

$\bullet$  if $\nu(q)\times\mbox{\boldmath $a$}\not=\mbox{\boldmath $0$}$, then
$\displaystyle\liminf_{\tau\rightarrow\infty}\tau^2e^{-2\tau R_D(p)}I(\tau)>0$,

$\bullet$  if $\nu(q)\times\mbox{\boldmath $a$}=\mbox{\boldmath $0$}$, then
$\displaystyle\liminf_{\tau\rightarrow\infty}\tau^3e^{-2\tau R_D(p)}I(\tau)>0$.

\noindent
From these together with (6.2) and Lemma 4.2 and replacing $\tau$ with $\tau\sqrt{\rho/\mu}$, we finally obtain:

$\bullet$  if $\nu(q)\times\mbox{\boldmath $a$}\not=\mbox{\boldmath $0$}$, then
$$\displaystyle
\liminf_{\tau\rightarrow\infty}e^{2\tau\sqrt{\rho/\mu}(R_1- R_D(p))}\Vert\mbox{\boldmath $w$}_{s0}\Vert_{L^2(B)}^2>0,
$$

$\bullet$  if $\nu(q)\times\mbox{\boldmath $a$}=\mbox{\boldmath $0$}$, then
$$\displaystyle
\liminf_{\tau\rightarrow\infty}\tau e^{2\tau\sqrt{\rho/\mu}(R_1- R_D(p))}\Vert\mbox{\boldmath $w$}_{s0}\Vert_{L^2(B)}^2>0.
$$

\noindent
This completes the proof of Lemma 4.3.

\section{Further problems and possible applications}

Everything is reduced to realizing or prescribing the desired input heat flux $f$ 
on the surface of the body, which has the form
$$\begin{array}{lll}
\displaystyle
f(x,t)=\nabla v(x,t)\cdot\nu(x), & x\in\partial\Omega, & 0<t<T,
\end{array}
$$
where
$$\begin{array}{lll}
\displaystyle
v(x,t)=\frac{1}{(\sqrt{2\pi t}\,)^3}
\int_{B_{R_2}(p)\setminus B_{R_1}(p)}e^{-\frac{\vert x-y\vert^2}{4t}}(R_2-\vert y-p\vert)(R_1-\vert y-p\vert)\,dy,
& x\in\Bbb R^3, & t>0.
\end{array}
$$
How to realize such heat flux will be the next technical problem.  A possible way is a combination of
an approximation of the Neumann-to-Dirichlet map (for the heat equation) and the principle of the superposition
like a phased array system.

The $B_{R_2}(p)\setminus B_{R_1}(p)$ has a shell-type geometry.  How about the case when $B_{R_2}(p)\setminus B_{R_1}(p)$
is replaced with a solid torus or ellipsoid?

Needless to say, the idea of using the shell-type initial data can be applied also to the inverse boundary value problem 
governed by the wave equation as considered in \cite{EIV, EV}.

Apply the idea to a dissipative medium as considered in \cite{E0}, viscoelastic medium
in \cite{IH}, etc.

$$\quad$$

\centerline{{\bf Acknowledgments}}

The author was partially supported by Grant-in-Aid for
Scientific Research (C)(No. 17K05331) and (B)(No. 18H01126) of Japan  Society for
the Promotion of Science. 
 
Some of this work was started by the author during his visit to University of Helsinki in March 2019.
The author would like to thank Samuli Siltanen for having useful discussions during the stay.

$$\quad$$

\vskip1cm
\noindent
e-mail address

ikehata@hiroshima-u.ac.jp

\end{document}